\documentclass[12pt]{article}
\usepackage{amsmath}%
\usepackage{amsfonts}%
\usepackage{amssymb}%
\usepackage{graphicx}
\usepackage[all]{xy}
\usepackage{amsthm}
\usepackage{latexsym}

\begin{document}

{\theoremstyle{plain}
  \newtheorem{theorem}{Theorem}[section]
  \newtheorem{corollary}[theorem]{Corollary}
  \newtheorem{proposition}[theorem]{Proposition}
  \newtheorem{lemma}[theorem]{Lemma}
  \newtheorem{question}[theorem]{Question}
  \newtheorem{conjecture}[theorem]{Conjecture}
}

{\theoremstyle{definition}
  \newtheorem{definition}[theorem]{Definition}
  \newtheorem{remark}[theorem]{Remark}
  \newtheorem{example}[theorem]{Example}

}

\numberwithin{equation}{section}

\bibliography{file 1, file 2, ....}
\bibliographystyle{style name}\bibliography{file 1, file 2, ....}

\bibliographystyle{style name}

\title{   On Unimodality of Hilbert Functions of Gorenstein Artin Algebras of Embedding Dimension Four  }

\author{Sumi Seo, Hema Srinivasan \\\small{Department of Mathematics, University of Missouri-Columbia} \\\small{E-mail:
sscz8@mail.missouri.edu,
srinivasanh@missouri.edu}}
\date{}

\maketitle

\begin{abstract}
We prove that the Hilbert functions of codimension four graded
Gorenstein Artin algebras $R/I$ are unimodal provided $I$ has a
minimal generator in degree less than five.  It is an open
question whether all Gorenstein $h$-vectors in codimension four are
  unimodal. In
this paper, we prove that Hilbert functions of all artinian
codimension four Gorenstein algebras starting with $(1,4,10, 20,
h_4,\cdots  )$, where $h_4\leq 34$ are unimodal.  Combining this
with the previously known results, we obtain that all Gorenstein
h-vectors $(1, h_1, h_2,h_3, h_4, \cdots) $ are unimodal if
$h_4\le 34$.      \end{abstract}


\newcommand\sfrac[2]{{#1/#2}}

\newcommand\cont{\operatorname{cont}}
\newcommand\diff{\operatorname{diff}}


\section{Introduction}%

 We study the  problem of whether the Hilbert function of a  standard graded Gorenstein algebra  is an SI sequence or at least whether it is unimodal.

In \cite{Stanley2}, Stanley proved that all Gorenstein $h$-vectors
are SI-sequences for $h_1 \leq 3$ using Buchsbaum-Eisenbud
structure theorems.  Bernstein and Iarrobino in \cite{B-I}
found examples of non-unimodal Gorenstein $h$-vectors for $h_1 \geq
5$, so it turns out that the unimodality  for Hilbert functions of Gorenstein algebras remained a question only in codimension four.
 In codimension four,  A. Iarrobino and
H. Srinivasan showed in \cite{I-S} that all Gorenstein h-vectors with
$h_2 \leq 7$ are SI-sequences and hence unimodal and J. Migliore, U. Nagel and F.
Zanello showed in \cite{M-N-Z} that Hilbert functions of all
artinian codimension four Gorenstein algebras starting with
$(1,4,h_2, h_3,\ldots)$, where $h_4\leq 33$ are SI-sequences.  It still remains open whether all
 Gorenstein h-vectors are SI-sequences, even  whether they must be
unimodal.  Since the Hilbert function $H_{R/I}$ is an O-sequence,
$h_4 \leq 35$. If $I$ has any minimal generator in degrees less
than five, $h_4\le 34$.  Thus, the only remaining case was $h_4= 34$.   We settle this
remaining case to prove:

\begin{theorem}
Let $R=k[x_1,x_2,x_3,x_4]$ where $k$ is a field of characteristic zero and $I$ be an artinian Gorenstein ideal of codimension four.
If a minimal generating set of $I$ contains  at least one generator of degree less than five,  then the $h$-vector of $R/I$ is
unimodal.
\end{theorem}

This is the main result of this paper. In section 2, we outline some of the preliminary results.
  In section 3, we show first that an $h$-vector of
an artinian Gorenstein algebra in codimension four with initial
degree 4 is unimodal  in theorem \ref{mainthm}, and  section 4 has some examples.

There has been steady work on the problem of
determining Gorenstein $h$-vectors and also whether they are unimodal beginning with Stanley's first counter example for nonunimodal $h$-vector in embedding dimension 13 followed by Bernstein and Iarrobino in embedding dimension 5 and higher.  It turns out that the examples of nonunimodal Gorenstein Hilbert functions all have $h_2$ less than the maximum predicted by Macaualay's bound for it.  So, the fact that any possibly nonunimodal Hilbert function in embedding dimension four has the maximum possible $h_2,h_3,h_4$ is interesting and gives hope for the conjecture that
Gorenstein Artin algebras of embedding dimension four are unimodal.

\section{Preliminaries}
Let $R=k[x_1, \cdots , x_n]$ be a polynimial ring over a field $k$
with standard grading and $k$ be an algebraically closed field of
characteristic zero. Let $I=\oplus_{i\geq 0}I_i$ be a homogeneous
ideal of $R$ and consider an artinian algebra $A=R/I=\oplus_{i\geq
0}A_i$. The Hilbert function of $A$ is $H_{A}:
\textbf{Z}{\longrightarrow} \textbf{Z}$ defined by
\begin{eqnarray*}
H_{A}(i)=dim_{k}R_i -dim_{k}I_i =\binom{n-1+i}{n-1} -dim_{k}I_i = h_i
\end{eqnarray*}

then   $h(A)=h=(h_0, h_1, h_2, \cdots, h_e)$ is known as the
$h$-vector of $A$ where $h_i=H_{A}(i)$ and $e$, the socle degree of
$A$, is the largest $i$ for which $dim_{k}A_{i}>0$. Without loss
of generality, we may assume that $I$ has no linear terms and
hence $h_1 =n$ is the embedding dimension.     If  $A$ is
Gorenstein, the $h$-vector is symmetric about the middle
$\left\lfloor \frac{e}{2} \right\rfloor$. An $h$-vector is said to
be $unimodal$ if it no longer increases once it starts decreasing.

If $h$ and $i$ are positive integers, then $h$ can be written uniquely
in the form of i-binomial expansion of $h$:
\begin{eqnarray*}
h=\binom{n_i}{i}+\binom{n_{i-1}}{i-1}+ \cdots +\binom{n_j}{j},
\end{eqnarray*}
where $h_i > h_{i-1}> \cdots > n_{j} \geq j \geq 1$.\\
We now define a collection of functions $-^{<i>}:
\textbf{Z}{\longrightarrow} \textbf{Z}$ as follows: if $h \in
\textbf{Z}$ has i-binomial expansion as above, then
\begin{eqnarray*}
h^{<i>}=\binom{n_i +1}{i+1}+\binom{n_{i-1}+1}{i}+ \cdots +\binom{n_j +1}{j+1}.
\end{eqnarray*}

A sequence of nonnegative integers ${\left\{h_i \right\}}_{i \geq
0}$ is called an {\it O-sequence} if $h_0 =1$ and $h_{i+1}\leq
h_{i}^{<i>}$ for all $i$.  Such sequences are precisely the ones
that occur as Hilbert functions of standard graded algebras  \cite
{Stanley1}. We say that maximal growth of the Hilbert function of
$A$ occurs in degree $i$ if $h_{i+1} = h_{i}^{<i>}$. The first
difference of the Hilbert function H  is denoted by $\Delta
H_{A}(i):=H_{A}(i)-H_{A}(i-1)=h_i-h_{i-1}$ for all $i \geq 1$. An
$h$-vector is called a {\it differentiable O-sequence up to $j$} if
$\Delta H_{\leq j}=(h_0, h_1 -h_0, h_2 -h_1, \cdots,
h_{j}-h_{j-1})$ is an O-sequence. An $h$-vector is called an {\it
SI-sequence(Stanley-Iarrobino sequence)} if it is symmetric about
the middle and the first half difference $\Delta H_{\leq
\frac{e}{2}}=(h_0, h_1 -h_0, h_2 -h_1, \cdots, h_{\left\lfloor
\frac{e}{2}\right\rfloor}-h_{{\left\lfloor \frac{e}{2}
\right\rfloor}-1})$ is an O-sequence.  The importance of this name
is due to the fact that in codimension three, being an SI sequence
is precisely the same as being Gorenstein.  However, it is not so
in codimension five and higher.  This leads us to concentrate on
codimension four.

We define $h_{<i>}$:
\begin{eqnarray*}
h_{<i>}=\binom{n_i -1}{i}+\binom{n_{i-1}-1}{i-1}+ \cdots +\binom{n_j -1}{j}.
\end{eqnarray*}

The following lemma encodes the change in  the Hilbert function
when one mods out by a sufficiently general linear form.  This is
used to cut down the dimension of $R$.

\begin{theorem}  \cite{Green} \label{Green}
If $l \in R$ is a general linear form, then the $h$-vector of $A/lA$
satisfies
\begin{center}
$h_{A/lA}(i)\leq (h_{A}(i))_{<i>}$.
\end{center}
\end{theorem}


If $R/I$ is a Gorenstein Artin algebra and $f \notin I$, then
$R/(I:f)$ is again Gorenstein, for if $I= ann F$ for some $F \in
k_{DP}[X_1, \cdots X_n]$ where $k_{DP}[X_1, \cdots X_n]$ is the
divided power algebra, then $(I:f) = ann F$.  Thus, $f$ is
homogeneous of degree $d$ and the socle degree of $R/(I:f)$ is
$e-d$. Thus,

\begin{lemma}\label{lem1}If $f \in R_d$ is a homogeneous form of degree $d$ and $f \notin
I$, then $R/(I:f)$ is Gorenstein with socle degree $e-d$, where $e$ is the socle degree of $R/I$.
\end{lemma}

 We apply this for a general linear form $l$ in $R/I$.  One of the first things we do is to compare the Hilbert function of $R/I$ with 
 $H_{R/(I:l)}$ for a suitably general $l$ and set up an induction on the socle degree.    Following observation, lemma \ref {lemm1}   shows how to compare the two $h$- vectors.  It is also proved in \cite{M-N-Z}. We
include a short proof.

 \begin{lemma}\label{lemm1}
Let $A=R/I$ be a graded Artin algebra with $dim_k A_1 \geq 2$.
Assume that $I_t$ has a GCD, $D$, of degree $d>0$. Set
$B=R/(I:D)$. Then
  If $i \leq t$, then 
  $$
  H_B (i-d)=H_A (i)- \left[
\binom{i+n-1}{n-1}- \binom{i-d+n-1}{n-1}\right].
$$
  Further,  if the
$h$-vector of B is non-decreasing (resp. increasing) in degrees
$\leq t-d$, then the $h$-vector of $A$ is non-decreasing (resp.
increasing) in degree $\leq t$.
\end{lemma}
\begin{proof}
From the exact sequence,
$$0\to R/(I:D)(-d)\to R/I \to R/(I,D) \to 0$$

we get $H_B (i-d)=H_{R/I} (i)- H_{R/(D)} (i)$ where
$(I,D)_i=(D)_i$ if $i\leq t$. So for $i\le t$,  $$H_B (i-d)=H_A
(i)- \left[ \binom{i+n-1}{n-1}- \binom{i-d+n-1}{n-1}\right]$$

Then, 
$$H_B(i-d+1) - H_B(i-d) = H_A(i+1)-H_A(i) -  \left[
\binom{i+n-1}{n-2}- \binom{i-d+n-1}{n-2}\right] $$ which proves the
lemma.
\end{proof}

From here on,  we let $R = k[x_1, x_2, x_3, x_4]$ and $a=min\{j\in \textbf{Z}| I_j \neq 0 \}$ denote the initial degree of $I$.  The following  lemma  relates to the situation  when there is exactly one minimal generator in any initial degree.   

\begin{lemma}\label{thm1}
Let $I$ be any ideal of codimension four and minimally generated
in degrees $\geq a$. Suppose $I$ has exactly one generator in
degree $a$ and all other generators in degree $\geq a+m$ for a
certain positive number $m$.  If $l_1$ and $l_2$ are two general
linear forms in $R$ and $(f_0,f_1,f_2, \cdots)$ is the Hilbert
function of $R/(I, l_1,l_2)$. Then $f_{a-1} =f_a = \cdots =
f_{a+m-1}=a$ and $f_{a+m}\leq a-1$. In particular, $\{ f_j \}$ is
non-increasing from $j=a-1$.
\end{lemma}
\begin{proof}
$f=(f_0,f_1,f_2, \cdots)$ is the Hilbert function of $R/(I, l_1,
l_2)$ that is an artinian quotient of $R/(l_1,l_2)\cong k[x,y]$.
Let $\bar{R}=R/(l_1,l_2)$ and $\bar{I}$ be the restriction to
$\bar{R}$. Then the h-vector of $\bar{R}/\bar{I}$ is $f$. Let $F
\in I$ be the minimal generator of degree $a$. We may assume that
a second minimal generator $G$ of $I$ comes in degree $a+m(m\geq
1)$. Since $I$ has no generator of degree less than $a$,
$f_{a-1}=dim_k (\bar{R}/\bar{I})_{a-1}=a$.
\begin{flalign*}
f_j &=dim_k (\bar{R}/\bar{I})_j
=dim_{k}{\bar{R}}_j-dim_{k} {\bar{I}}_j =(j+1)-dim_{k}
{(\bar{F})}_j \\
&=(j+1)-(j-a+1) =a \ \ \text{for \ all}\  j, \ a \leq j < a+m.
\end{flalign*}
So, $f_{a-1} = \cdots = f_{a+m-1}=a$. The second minimal generator
$G$ is of degree $a+m$, then $(F, G)_{a+m}\subset I_{a+m}$.
$dim_{k} \bar{I}_{a+m} >
dim_{k}(\bar{F})_{a+m}=m+1$.
\begin{flalign*}
f_{a+m} &=dim_{k} \bar{R}_{a+t} - dim_{k}\bar{I}_{a+t} <
(a+m+1)-(m+1)=a=f_{a+m-1}
\end{flalign*}

For all $j\geq a-1$, $f_j \geq f_{j+1}$ because  $dim_{k}{\bar{I}}_j < dim_{k}{\bar{I}}_{j+1}$ unless $\bar{I}_{j+1}=0$. So, $\{ f_j \}$ is non-increasing from $j\geq a-1$.
\end{proof}

We note that if the Hilbert function of $R/I$ is $(1,4,10, 20,
h_4, \cdots )$, then the ideal $I$ is generated in degrees $4$ and
higher.  So, $a=4$.
 If  $m$ is the difference between the first and second minimal degrees of generators in $I$, we have $4$ entries in the $h$-vector $(f_0,f_1,f_2, \cdots)$ in degree 3 up to degree one less than the second minimal degree,
that is, $f_3 =f_4 = \cdots = f_{m+3}=4$.

The following lemma plays a very important role.  The example \ref{ex1} is one with this situation. 
\begin{lemma}\label{lemm2}
Let $J=(f,g,h)\subset I$ with $\deg h=\deg g +1$.
Suppose $D=GCD(f,g)$ is of degree $d(\deg f-2\leq d < \deg f)$ and $h \notin
(D)$. Let $l_1$ and $l_2$ be general linear forms. Then
      \begin{center}
       $H_{R/(J,l_1,l_2)}(j)=H_{R/(J,D,l_1,l_2)}(j)$ for $j\geq \deg h$.
     \end{center}
     \end{lemma}
\begin{proof}
Let $\bar{R}=R/(l_1 ,l_2)$.    $\bar{I}$ and $\bar{J} $  are the image of $I$ and $J$ respectively in $\bar{R}$. Since $D$
is the GCD of $f$ and $g$, $f=f'D$, $g=g'D$ and $(f',g')$ is the
regular sequence with $\deg f'=\deg f-d$ and $\deg g'=\deg g-d$.
Consider the exact sequence:
\[
0\,\rightarrow\,[\bar{R}/(\bar{J}:\bar{D})](-d)\, \rightarrow\,\
\bar{R}/\bar{J} \,\ \rightarrow \,\ \bar{R}/(\bar{J},\bar{D})\,\
\rightarrow\,\ 0.
\]

Since $(\bar{J}:\bar{D})=(\bar{f'},\bar{g'},\bar{h})$, starting
$j={\deg h}$
\begin{flalign}\label{formula1}
dim_k [\bar{R}/\bar{J}]_{\deg h} = dim_k
[\bar{R}/(\bar{f'},\bar{g'},\bar{h})]_{{\deg h}-d}+dim_k
[\bar{R}/(\bar{J},\bar{D})]_{\deg h}
\end{flalign}

$(\bar{f'},\bar{g'},\bar{h})$ has at most two minimal generators
up to degree ${\deg h}-d$. So we get
\begin{flalign*}
dim_k [\bar{R}/(\bar{f'},\bar{g'},\bar{h})]_{{\deg h}-d} &= dim_k
[\bar{R}/(\bar{f'},\bar{g'})]_{{\deg h}-d}\\ &= ({\deg h}-d+1)-({\deg h}- {\deg f}+1)-({\deg h}-{\deg g}+1)\\
&={\deg f}-d-2 \\
dim_k [\bar{R}/\bar{J}]_h &- dim_k [\bar{R}/(\bar{J},\bar{D})]_{\deg h}
={\deg f}-d-2 \leq 0.
\end{flalign*}

Thus, 
$$
dim_k [\bar{R}/\bar{J}]_{\deg h} \leq dim_k
[\bar{R}/(\bar{J},\bar{D})]_{\deg h}.
$$
 On the other hand, 
 $$
 (\bar{J},
\bar{D})=(\bar{D},\bar{f},\bar{g},\bar{h})=(\bar{D},\bar{h})\supset
(\bar{f},\bar{g},\bar{h})= \bar {J}.
$$
 It follows 
 $$
 dim_k
[\bar{R}/\bar{J}]_{\deg h} \geq dim_k [\bar{R}/(\bar{J},\bar{D})]_{\deg h}.
$$

If $j >\deg h$, $(I,D,l_1,l_2)_{j}=(I,l_1,l_2)_{j}$, which implies the two $h$-vectors are equal.\\
\end{proof}


\section{$H=(1,4,10,20,34, \cdots ,34,20,10,4,1)$  }
In this section,  $R= k[x_1, x_2, x_3, x_4]$ and  $I$ is  a grade
four Gorenstein ideal in $R$.  We consider the case $H_{R/I}  =
(1,4,10,20,  \cdots)$. So, $I$ is minimally generated in degree
$4$ and all other generators in higher degrees.

Let $l_1$ and $l_2$ be sufficiently general linear forms so that
$R/(I,l_1,l_2)$ is isomorphic to $k[x_1,x_2]/I$.  Recall that
$(I:l_1)$ is Gorenstein and $H_{R/I}$ denotes the Hilbert function
of $R/I$. From the exact sequence
 $$0\to R/(I:l_1) (-1)\to R/I \to R/(I,l_1) \to 0$$ we get,
$$H_{R/(I,l_1)}(i)  = H_{R/I}(i) - H_{R/(I:l_1)}(i-1),$$ and
$$H_{R/(I,l_1,l_2)}(i) = H_{R/(I,l_1)}(i) - H_{R/((I,l_1):l_2)}(i-1).$$

Let $b_i :=H_{R/(I:l_1)}(i-1)$, the $(i-1)$-th entry of the
Hilbert function of $H_{R/(I:l_1)}$, then we have
$H_{R/(I:l_1)}(-1)=(0,1,b_2,b_3, \cdots )$ and
$H_{R/(I,l_1,l_2)}=(f_0,f_1,f_2, \cdots)$. Thus, 
$$
f_i =
h_i-b_i-H_{R/((I,l_1):l_2)}(i-1) \le
(h_i)_{<i>}-H_{R/((I,l_1):l_2)}(i-1).
$$








\begin{theorem}\label{mainthm}
Let $R=k[x_1,x_2,x_3,x_4]$ and $I$ be an artinian Gorenstein ideal of codimension four.
If $I$ contains exactly one minimal generator of degree four, then the $h$-vector of $R/I$ is unimodal.
\end{theorem}
\begin{proof}
 Since $R/I$ has no generator in degrees 1, 2 and 3, its Hilbert function is given by $h=(1,4,10,20,h_4 \cdots,
h_e)$  with $h_4 = 34$.     $h_4=34$ because  there is exactly one minimal generator for $I$ in degree four.  We call it $f$.  Let the next minimal generator be in degree $4+m$.

 We use induction on the socle degree  of $R/I$. It is clearly unimodal when  the socle degree is 1 and in fact when it is less than or equal to $8$.
 Suppose all such Gorenstein Artin algebras with socle degrees less than $e$ have unimodal Hilbert functions.
 So, we may assume that the $h$-vector of $R/I$ has the least socle degree $e>8$.   Suppose the $h$-vector of $R/I$ is not unimodal.

 $R/(I:l_1)$ has socle degree $e-1$  and hence by induction $H_{R/(I:l_1)} = (b_1, b_2, \cdots b_e)$ is unimodal.
 By lemma \ref{thm1}, the $H_{R/(I,l_1,l_2)}$ is $(1,2,3,4,4,\cdots, 4, f_{4+m},\cdots)$ and non-increasing after $f_3$.

We note that
\begin{flalign*}
H_{R/(I,l_1)} &= H_{R/I}-H_{R/(I:l_1)}(-1) \\
              &=(1, 4, 10, 24, 34, \cdots, 34, 20, 10, 4, 1) - (0, 1, b_2, b_3, \cdots b_3, b_2, 1)\\
              &=(1, 3, 6,  10, 14, \cdots)\\ 
H_{R/(I,l_1, l_2)} &= H_{R/(I,l_1)} - H_{R/(I,l_1):l_2)}(-1)\\
H_{R/(I,l_1, l_2)} &= H_{R/I}-H_{R/(I:l_1) (-1)} - H_{R/((I,l_1):l_2)}(-1).
\end{flalign*}
So,
$$(1,2,3,4,4,\cdots, 4, f_{4+m},\cdots)=(1, 3, 6,  10, 14, \cdots) - H_{R/((I,l_1):l_2)}(-1)$$

Further, we have 
$$
b_t =H_{R/(I:l_1)}(t-1) \le H_{R/I}(t-1) =
h_{t-1}
$$
 and
 $$
 H_{R/((I,l_1):l_2)}(t-1)\le H_{R/(I,l_1)}(t-1) = h_{t-1}-b_{t-1}.
 $$

If $h=(h_0, h_1,h_2, \cdots, h_e)$ is not unimodal, there is the least integer $i \leq \frac{e}{2} -1$ such that $h_i
> h_{i+1}$ equivalent to $h_{e-1-i}<h_{e-i}$ by symmetry. We must have $i+1 \geq 5$, then $\frac{e}{2}\geq i+1 \geq 5$, so $e \geq 10$.
 Since $h_i >h_{i+1}$ and $h_{e-1-i}<h_{e-i}$, and $(1,b_2, b_3, \cdots) $ is unimodal, we get 
 $$
 f_{e-i} = h_{e-i}-b_{e-i}-H_{R/((I,l_1):l_2)}(e-i-1)>h_{e-i-1}-b_{e-i-1} -H_{R/((I,l_1):l_2)}(e-i-1)\ge 0.
 $$
 Recall that $\frac{e}{2}+1 \leq e-i$.  Thus, we have two cases: $\frac{e}{2}+1 < e-i$  and $\frac{e}{2}+1 = e-i$.
\vskip .2truein
\noindent{\bf Case I}:\ $\frac{e}{2}+1 < e-i$

 Since $e-i-2 > \frac{e}{2}-1 \geq 4$, definitely $f_{e-i-2}\leq 4$.  $f_{e-i}>0$ allows the possible values of the tuple $(f_{e-i-2},f_{e-i-1},f_{e-i})$ as follows:
\begin{itemize}
\item (4,4,4),(4,4,3),(4,4,2),(4,4,1),(4,3,3),(4,2,2),(4,1,1),\\
      (3,3,3),(3,3,2),(3,3,1),(3,2,2),(3,1,1),\\
      (2,2,2),(2,2,1),(2,1,1),(1,1,1)\\
      \vskip .1truein
      We know that $f_{e-i-2}=1,2,3$ or 4 and $e-i-2>4$, so $f_{e-i-2} \leq e-i-2$. When $f_{e-i-2}=f_{e-i-1}$,
      $f_{e-i-2}^{\ \ \ \ <e-i-2>}=f_{e-i-2}=f_{e-i-1}$ which implies that there is maximal growth of Hilbert function of $R/(I,l_1,l_2)$in degree $e-i-2$.
      When $f_{e-i-1}=f_{e-i}=1, 2, 3$ or 4, $e-i-1 \geq 4$ and $f_{e-i-1}^{\ \ \ \ <e-i-1>}= f_{e-i}=f_{e-i-1}$ allow us to have
      maximal growth in Hilbert function in degree $e-i-1$. So $\bar{I}_{e-i-1}$ has a GCD of appropriate degree by proposition 2.7 in \cite{B-G-M}.
      By lifting to $I$, $I_{e-i-1}$ has a GCD. For, if $I_t$ has no GCD then there are $a, b$ in $I_t$ forming a regular sequence. They are still regular modulo general enough $l_1$, $l_2$, so $\bar{I}_{t} = I_t$ in $R/(l_1,l_2)$ has no GCD(see \cite{Davis}). Further codimension of $I_t$ = codimension of $\bar{I}_t$ =1. Hence $I_{e-i-1}$ has a GCD of degree 1,2,3 or 4. By $e-i-1> \frac{e}{2}$ and lemma \ref{lemm1}, the $h$-vector of $R/I$ is unimodal, contradicting our assumption of nonunimodality.
\item (4,3,2),(4,3,1),(4,2,1)\\
      By theorem \ref{thm1}, the second minimal generator of $I$ comes in degree $e-i-1(> \frac{e}{2})$, so
      the unimodality follows immediately.
\item (3,2,1)\\
      $f_{e-i-3}$ is either 4 or 3. If $(f_{e-i-3},f_{e-i-2},f_{e-i-1},f_{e-i})=(4,3,2,1)$, then the second minimal
      generator comes in degree $e-i-2(>\frac{e}{2}-1)$. Actually, $e-i-2 \geq
      \frac{e}{2}$. So the $h$-vector of $R/I$ is non-decreasing to degree $\frac{e}{2}$, leading to contradiction.
      If $(f_{e-i-3},f_{e-i-2},f_{e-i-1},f_{e-i})=(3,3,2,1)$, then $I_{e-i-2}$ has a
      GCD, D, of degree 3. Note that by induction hypothesis, $R/(I:D)$ is unimodal and has an
      $h$-vector that is non-decreasing up to degree $\frac{e-3}{2}$ since the socle
      degree of $R/(I:D)= e-3$. So it is non-decreasing up to $\frac{e}{2}-3(\leq (e-i-2)-3)$.
      Hence, lemma \ref{lemm1} gives that the h-vector of $R/I$ is non-decreasing up to
      degree $\frac{e}{2}$, and is unimodal by symmetry.
\end{itemize}

\noindent {\bf Case II}:\ $\frac{e}{2}+1 = e-i$

This case is of $i=\frac{e}{2}-1$ and $e$ is even. Consequently,
the non-unimodality occurs in the very middle, that is
$h_{\frac{e}{2}-1} > h_{\frac{e}{2}} < h_{\frac{e}{2}+1}$.
Consider the tuple
$(f_{\frac{e}{2}-1},f_{\frac{e}{2}},f_{\frac{e}{2}+1})$ of $f$.
With the same argument above, we get rid of the tuples where two
entries are equal since they have a GCD in degree $\frac{e}{2}$
and apply lemma \ref{lemm1}. The cases that we need to deal
with are (4,3,2),(4,3,1),(4,2,1) and (3,2,1).

\begin{itemize}
\item (4,3,2),(4,3,1),(4,2,1)\\
      By theorem \ref{thm1}, the second minimal generator of $I$ comes in degree $\frac{e}{2}$, so
      the unimodality follows immediately.
\item (3,2,1)\\
     We claim that $I_{\frac{e}{2}}$ has no GCD. Suppose not, $I_{\frac{e}{2}}$ has a GCD, D, of
      degree $d (0<d<3)$. $(I,D)_j =(D)_j$ if $j \leq \frac{e}{2}$.
\begin{flalign*}
\Delta h_{R/(I,D)}(\frac{e}{2})= \Delta h_{R/(D)}(\frac{e}{2})= \binom {2+ {\frac{e}{2}}}{2} - \binom {2+ {\frac{e}{2}}-d}{2} \geq \frac{e}{2} +1
\end{flalign*}
Note that $H_{R/I}(i)=H_{R/(I:D)}(i-d)+H_{R/D}(i)$ for all $i>d$. Since $h_{\frac{e}{2}} < h_{\frac{e}{2}-1}$,
\begin{flalign*}
0&>h_{\frac{e}{2}}-h_{\frac{e}{2}-1}\\
 &= H_{R/(I:f)}(\frac{e}{2}-d) - H_{R/(I:f)}(\frac{e}{2}-1-d)+ \Delta h_{R/(f)}(\frac{e}{2})\\
 &\geq H_{R/(I:f)}(\frac{e}{2}-d) - H_{R/(I:f)}(\frac{e}{2}-1-d) + \frac{e}{2} +1\\
&\Rightarrow  \  H_{R/(I:f)}(\frac{e}{2}-d) < H_{R/(I:f)}(\frac{e}{2}-1-d)
\end{flalign*}
This implies contradiction to $H_{R/(I:f)}(\frac{e}{2}-d) \geq H_{R/(I:f)}(\frac{e}{2}-1-d)$ by unimodality of Gorenstein $R/(I:f)$. This completes the proof of this claim.

      $f_{\frac{e}{2}-2}$ is either 4 or 3. If $(f_{\frac{e}{2}-2},f_{\frac{e}{2}-1},f_{\frac{e}{2}},f_{\frac{e}{2}+1})=(3,3,2,1)$,
      then $I_{\frac{e}{2}-1}$ has a GCD, A, of degree 3. The second minimal generator of $I$ comes in degree $4+m$ for some $m(0<m<\frac{e}{2}-4)$. In fact, $I$ has exactly one minimal
      generator of degree $4+m$. Otherwise, $f_{4+m}< 3$, not satisfying this case  $(f_{\frac{e}{2}-1},f_{\frac{e}{2}},f_{\frac{e}{2}+1})=(3,2,1)$.
      Denote $g$ the second minimal generator of $I$ then $f,g \in (A)$, so we get $(I,A)_j =(A)_j$ if $j \leq \frac{e}{2}-1$.
      Since the $(f_0,f_1,f_2, \cdots)$ has dropped at degree $\frac{e}{2}$ by one, the third minimal generator $h$ comes in degree $\frac{e}{2}$. Since $I_{\frac{e}{2}}$ has no GCD, $h \notin (A)$.
      $(I,A)=(A,f,g,h, \cdots)=(A,h, \cdots)$ implies $(I,A)$ has at most two minimal generators up to degree $\frac{e}{2}$ and its second minimal generator is $h$ of degree $\frac{e}{2}$. By lemma \ref{thm1}, the following table is obtained:
      \vskip .1truein

\begin{small}
\begin{tabular}{c|c c c c c c c c c c c c c c }
      Degree                   &  0 & 1 & 2  & 3  & 4  & 5      & $\cdots$ & $3+m$ & $4+m$ & $\cdots$  & $ \frac{e}{2}  -1$ & $ \frac{e}{2}$\\ \hline
 $H_{R/(I, l_1, l_2)}$    &  1 & 2 & 3  & 4  & 4  & 4      & $\cdots$ & 4     &3      &  $\cdots$ &                                           3 & 2 \\ \hline
 $H_{R/(A,I,l_1,l_2)}$    &  1 & 2 & 3  & 3  & 3  & 3      & $\cdots$ & 3     &3      &  $\cdots$ &                                           3 & 2 \\
\end{tabular}
\end{small}
            \begin{flalign*}
      H_{R/(I,l_1,l_2)}(j)&=H_{R/(I,A,l_1,l_2)}(j)=3 \  \text{for} \ 4+m < j < \frac{e}{2}-1\\
      H_{R/(I,l_1,l_2)}(\frac{e}{2})&=H_{R/(I,A,l_1,l_2)}(\frac{e}{2})=2,
      \end{flalign*}

     which enables us to get $(I,A,l_1,l_2)_{\frac{e}{2}}=(I,l_1,l_2)_{\frac{e}{2}}$. We apply the information obtained from lifting  $(I,A,l_1,l_2)/(l_1,l_2)$ to
      $(I,A)$ as we have done before. The Hilbert function of $R/(I,A)$ is non decreasing up to degree $\frac{e}{2}$, for $(I,A)$ has at most two minimal generators up to degree $\frac{e}{2}$.
      Note that by induction hypothesis, $R/(I:A)$ is unimodal and has an
      h-vector that is non-decreasing up to degree $\frac{e-3}{2}$ since the socle
      degree of $R/(I:A)= e-3$. So it is non-decreasing up to $\frac{e}{2}-3$.
      Hence, the h-vector of $R/I$ is non-decreasing up to degree $\frac{e}{2}$ by the relationship between the Hilbert functions. So, it is unimodal by symmetry.

     If $(f_{\frac{e}{2}-2},f_{\frac{e}{2}-1},f_{\frac{e}{2}},f_{\frac{e}{2}+1})=(4,3,2,1)$, then the second minimal
     generator comes in degree $\frac{e}{2}-1$. Let us denote G the second minimal generator of degree $\frac{e}{2}-1$,
     then $I=(f, g, \text{others})$. If there is no generator of degree $\frac{e}{2}$ in
     $I$, $I$ has at most two generators up to degree $\frac{e}{2}$
     and it is obviously unimodal. We may assume that $I$ has exactly one generator of degree $\frac{e}{2}$ denoted by $h$. If $(f, g)$ forms a regular sequence,
     there is no generator of degree $\frac{e}{2}$ with the condition $(f_{\frac{e}{2}-1},f_{\frac{e}{2}},f_{\frac{e}{2}+1})=(3,2,1)$. So the h-vector is unimodal. We consider that $f, g$ are not a regular sequence in which there is a GCD, D, of $f, g$ of degree $d (< degf).$

     Let $J=(f,g,h)\subset I$, then $(J,D)=(D,h)$ has the second generator $H$ of degree $\frac{e}{2}$, so the h-vector $\bar{R}/(\bar{J},\bar{D})$ has dropped
     at degree $\frac{e}{2}$ by lemma \ref{thm1}. We have the h-vectors of $\bar{R}/(\bar{J},\bar{D})$ as
     follows:
\vskip .1truein
\begin{center}
\begin{tabular}{c|c|c c c c c c c c c c c c c c }
                                &           & 0 & 1 & 2  & 3  & 4  & 5     & $\cdots$ & $ \frac{e}{2}-2$ & $\frac{e}{2} -1$ & $\frac{e}{2} $\\ \hline
 $H_{\bar{R}/\bar{J}}$              &           & 1 & 2 & 3  & 4  & 4  & 4     & $\cdots$ &  4 &  3 & 2 \\ \hline
 $H_{\bar{R}/(\bar{J},\bar{D})}$    & d = 3 & 1 & 2 & 3  & 3  & 3  & 3     & $\cdots$ &  3 &  3 & 2 \\
 $H_{\bar{R}/(\bar{J},\bar{D})}$    & d = 2 & 1 & 2 & 2  & 2  & 2  & 2     & $\cdots$ &  2 &  2 & 1 \\
 $H_{\bar{R}/(\bar{J},\bar{D})}$    & d = 1 & 1 & 1 & 1  & 1  & 1  & 1     & $\cdots$ &  1 &  1 & 0 \\
\end{tabular}
\end{center}
\vskip .1truein
      When $d=2,3$,  we apply lemma \ref{lemm2} to $\deg f=4$, $\deg g=\frac{e}{2}-1$ and $\deg h=\frac{e}{2}$. We get, 
       $$H_{R/(J,l_1,l_2)}(\frac{e}{2})=H_{R/(J,D,L_1,L_2)}(\frac{e}{2})=2.$$
     
     Since $J_{\frac{e}{2}}=I_{\frac{e}{2}}$, we conclude that
       $$(I,l_1,l_2)_{\frac{e}{2}}=(I,D,l_1,l_2)_{\frac{e}{2}}.$$     

     If $d=3$, we apply the information of lifting from $(I,D,l_1,l_2)/(l_1,l_2)$ to
      $(I,D)$. It is unimodal since the proof remaining is similar to $(I,A)$  in the
 case of $(f_{\frac{e}{2}-2},f_{\frac{e}{2}-1},f_{\frac{e}{2}},f_{\frac{e}{2}+1})=(3,3,2,1)$.

     If $d=2$, from table above $H_{\bar{R}/\bar{J}}(\frac{e}{2}
)=2 > 1=H_{\bar{R}/(\bar{J},\bar{D})}(\frac{e}{2})$, which
contradicts the result of lemma \ref{lemm2}. The case of $d=2$ can
not occur.

     If $d=1$,  in $R/(J,D,l_1,l_2)\cong k[x]/\bar{h}$
     where $\deg \bar{h}=\frac{e}{2}$. So by theorem \ref{thm1},
     $H_{R/(J,D,l_1,l_2)}(j)=0$ for $j \geq \frac{e}{2}$ and $H_{\bar{R}/(\bar{J},\bar{D})}(\frac{e}{2})=0$. Note that $f=Df'$ and $g=Dg'$ where $(f',g')$ forms a regular sequence with $\deg f'=3$ and $\deg g'=\frac{e}{2}-2$.
     In (\ref{formula1}),
\begin{flalign*}
H_{\bar{R}/(\bar{J}:\bar{D})}(\frac{e}{2}-1) &=
H_{\bar{R}/\bar{J}}(\frac{e}{2}) - H_{\bar{R}/(\bar{J},\bar{D})}(\frac{e}{2}) \\
&=2-0 = 2
\end{flalign*}

However, since $(f',g')$ forms a regular sequence we get the contradiction as follows:
\begin{flalign*}
H_{\bar{R}/(\bar{J}:\bar{D})}(\frac{e}{2}-1)&=H_{\bar{R}/(\bar{f'},\bar{g'},\bar{h})}(\frac{e}{2}-1)\\
&=H_{\bar{R}/(\bar{f'},\bar{g'})}(\frac{e}{2}-1)=1
\end{flalign*}

Hence the condition
$(f_{\frac{e}{2}-2},f_{\frac{e}{2}-1},f_{\frac{e}{2}},f_{\frac{e}{2}+1})=(4,3,2,1)$
and the case
of $d=1$ do not occur simultaneously. We complete the proof.
\end{itemize}
\end{proof}

 \begin{theorem} Let $R=k[x_1,x_2,x_3,x_4]$ where $k$ is a field of characteristic zero and $I$ be an artinian Gorenstein ideal of codimension four.   If $I$ has at least one generator in degree less than five, then the $h$-vector of $R/I$ is unimodal.
 \end{theorem}
\begin{proof}
 Since $R/I$ has  at least one generator in degree less than five, the $H_{R/I}  =(1,h_1,h_2, h_3 ,h_4,  \cdots,
h_e)$  with $h_4 \leq 34$. If $h_4 \leq 33$,
it is unimodal \cite{M-N-Z}. If $h_4= 34$, then 
$H_{R/I}  =(1,4,10, 20, 34, \cdots,
h_e)$.  Then we have just proved in the above theorem \ref{mainthm}  that $H_{R/I}$ is unimodal.
\end{proof}
\section{Examples}
Here are some examples of Gorenstein Artin algebras.  
\begin{example}\label{ex1}

  $$I=(x^2 w^2, x^6, x^4 y^3  - zw^6, w^8, y^9 w^2, z^{11}, x^2 z^{10}, y^{12}w, y^{12}z, y^{13}-x^3 z^9 w,x^3 y^{12}, y^9 z^{10})$$

 Its $h$-vector is 

   $(1, 4, 10, 20, 34, 52, 73, 95, 116, 136, 156, 174, 187, 191, 187, 174, \cdots ,52, 34, 20, 10, 4, 1). $ 
 Since  $h_4=34$ and $h_5 \leq 52$,  there is maximal growth in $H_{R/I}$  in degree 4. If we take,$l_1=2x-5y+13z-7w$ and $l_2=-11x-4y+5z+9w$, then the Hilbert function of $R/(I,l_1,l_2)$ is $$H_{R/(I,l_1,l_2)}= (1, 2, 3, 4, 4, 4, 3, 1, 0, 0, \cdots).$$ 
 
 Note that $D=x^2$ is the gcd of the first two generators and the third generator is not divisible by $x^2$.   Further, the third generator is of degree $7$, one more than the degree of the second generator.  This is the situation in the lemma \ref{lemm2}  
 
 We have $$  H_{R/(J,l_1,l_2)}= (1, 2, 3, 4, 4, 4, 3, 1, 0, 0, 0, \cdots).$$ 
 The Hilbert function of $R/(J,D,l_1,l_2)$ is $(1, 2, 2, 2, 2, 2, 2, 1, 0, 0, 0, \cdots)$.  
\end{example}

\begin{example}
Let  ideal $I$ be generated by 
$$(y^2 w^2 , y^4 w, y^4 z, xy^4 , w^8 , x^6 y^2 , z^{10}w^2 , y^2 z^{11}   - x^7 w^6 , x^{13}, x^6 z^{10}, z^{21}, y^{26}- x^5 z^{20}w).$$

The $h$-vector for this $R/I$ is 

     $(1, 4, 10, 20, 34, 49, 66, 85, 104, 121, 137, 153, 168, 179, 184, \\\indent  184, 179,   168, 153, 137, 121, 104, 85, 66, 49, 34, 20,10, 4, 1). $   

If we take  $l_1=x - 3y + 15z - 2w$ and $l_2=- 13x - 4y+ 5z + 8w$  as general enough linear forms, then the $h$-vector of $R/(I,l_1,l_2) $ is 
  $(1, 2, 3, 4, 4, 2, 2, 2, 1 ).$      
The maximal growth in Hilbert function of $R/(I,l_1,l_2)$ occurs in degree 5 and 6. Then $I_5, I_6$ and $I_7$ have a gcd of degree 2.  
\end{example}

The next two examples are some Gorenstein ideals generated in degrees 5 and higher and they 
do have unimodal Hilbert functions.  

\begin{example}
Let  ideal $I$ be generated by 
$$(x^2 w^4 , x^6, z^3w^4, xy^6,z^8,w^9, y^6z^3,x^5z^5-y^5w^5, y^{11}) .$$

The $h$-vector for this $R/I$ is 
  $$(1, 4, 10, 20, 35, 55, 79, 104, 127, 137, 143, 149,143,137, 127,104,79,55,35,20,10,4,1 ).$$   

 \end{example}
 \begin{example}
Let  ideal $I$ be generated by 
$$(y^3 w^3,y^5 w,y^5z,xy^5,w^9,x^7y^3,y^2z^{11}-x^7w^6,z^{11}w^3, x^{14}, x^7z^{11},z^{22},y^{29}-x^6z^{21}w^2) .$$

The $h$-vector for this $R/I$ is 

   $(1, 4, 10, 20, 35, 56, 80, 107,   137, 169 , 201,231,259, 285,307,322,329,\\\indent 329,322,307,\cdots   \cdots  ,56,35,20,10,4,1 ). $   

 \end{example}


\end{document}